\documentclass{amsart}
\usepackage{amsmath}
\usepackage{amssymb}
\usepackage{amsthm}
\usepackage{amscd}
\usepackage{amsfonts}
\usepackage{graphicx}%
\newcommand{\set}[1]{\left\{#1\right\}}

\theoremstyle{plain}
\newtheorem{thm}{Theorem}[section]

\newtheorem{lem}[thm]{Lemma}
\newtheorem{prop}[thm]{Proposition}

\theoremstyle{definition}

\def\tr{\operatorname{tr}}

\begin{document}
\title[Differential Recursion Relations for Laguerre Functions]{Differential Recursion Relations for Laguerre Functions on Hermitian Matrices}
\author{Mark Davidson and Gestur \'{O}lafsson}
\address{Department of Mathematics, Louisiana State University, Baton Rouge, LA\ 70803, USA
}
\email{davidson@math.lsu.edu}
\email{olafsson@math.lsu.edu}
\thanks{Research by G. \'{O}lafsson supported by NSF grant DMS-0070607 and DMS-0139783}

\begin{abstract}
In our previous papers  \cite{doz1,doz2}  we studied Laguerre
functions and polynomials on symmetric cones $\Omega=H/L$. The
Laguerre functions $\ell^{\nu}_{\mathbf{n}}$,
$\mathbf{n}\in\mathbf{\Lambda}$, form an orthogonal basis in
$L^{2}(\Omega,d\mu_{\nu })^{L}$ and are related via the Laplace
transform to an orthogonal set in the representation space of  a
highest weight representations $(\pi_{\nu}, \mathcal{H}_{\nu})$ of
the automorphism group $G$ corresponding to a  tube domain
$T(\Omega)$.  In this article we consider the case where $\Omega$
is the space of positive definite Hermitian matrices and
$G=\mathrm{SU}(n,n)$. We describe the Lie algebraic realization of
$\pi_{\nu}$ acting in $L^{2}(\Omega,d\mu_{\nu})$ and use that to
determine explicit differential equations and recurrence relations
for the Laguerre functions.

\end{abstract}
\maketitle

\section*{Introduction}

In this article we continue our study of Laguerre functions and Laguerre
polynomials on symmetric cones. Let $\nu>0$ and $\alpha=\nu-1$. It is an old
and well known fact that the classical Laguerre functions $\ell_{n}^{\alpha
}(t)=e^{-t}L_{n}^{\alpha}(2t)$ form an orthogonal basis for $L^{2}%
(\mathbb{R}^{+},t^{\alpha}dt)$ and satisfy several differential equations and
recursion relations. The same holds for the Laguerre polynomials
$L_{n}^{\alpha}(t)$, in fact they may be defined as the polynomial solution to
a second order differential operator known as Laguerre's differential
operator. In our previous article \cite{doz1} it was shown that the
representation theory of $G=\widetilde{\mathrm{SL}(2,\mathbb{R})}$, the
universal covering group of $\mathrm{SL}(2,\mathbb{R})$, unified many of the
disparate properties of Laguerre functions. Briefly, the Laguerre functions
$\ell_{n}^{\alpha}$ form a basis of the space of $K=\widetilde{\mathrm{SO}%
(2)}$-finite vectors of a highest weight representations
$T_\alpha$ of $G$ acting on $L^{2}(\mathbb{R}^{+},t^{\alpha}dt)$.
The Lie algebra $\mathfrak{sl}(2,\mathbb{C})$ acts by at most
second order differential operators on the space of smooth vectors
of those representations. In
particular, the elements%
\[
X^{-}=\left(
\begin{array}
[c]{cc}%
1 & i\\
i & -1
\end{array}
\right)  ~,\quad X^{+}=\left(
\begin{array}
[c]{cc}%
-1 & i\\
i & 1
\end{array}
\right)  =-\overline{X^{-}},\quad\mathrm{and}\quad X^{0}=\left(
\begin{array}
[c]{cc}%
0 & -i\\
i & 0
\end{array}
\right)
\]
act, respectively,  by the differential operators%
\begin{align}
D^{-}  &  =tD^{2}+(2t+\nu)D+(t+\nu)\label{slaction1}\\
D^{+} &=tD^{2}-(2t-\nu)D+(t-\nu) \label{slaction2}\\
D^{0}  & =tD^{2}+\nu D-t.\label{slaction3}
\end{align}

Notice that $X^{0}$ is a basis for $\mathfrak{so}(2)$ and that our
normalization is such that%
\[
\lbrack X^{0},X^{-}]=2X^{-},\quad\lbrack X^{0},X^{+}]=-2X^{+},\quad
\mathrm{and}\quad\lbrack X^{+},X^{-}]=4X^{0}~.
\]
Furthermore, the action of $D^{0}$ gives us the following
differential equation for the Laguerre functions:
\begin{equation}
(tD^{2}+\nu D-t)\ell_{n}^{\alpha}=-(2n+\nu)\ell_{n}^{\alpha} \label{eq0intro}%
\end{equation}
and the actions of the annihilating operator $D^{-}$ and the
creation operator $D^{+}$ give us the recurrence relations for the
Laguerre functions:
\begin{equation}
D^{-}\ell_{n}^{\alpha}=(tD^{2}+(2t+\nu)D+(t+\nu))\ell_{n}^{\alpha}%
=-2(n+\nu-1)\ell_{n-1}^{\alpha} \label{eq1intro}%
\end{equation}
and
\begin{equation}
D^{+}\ell_{n}^{\alpha}=(tD^{2}-(2t-\nu)D+(t-\nu))\ell_{n}^{\alpha}%
=-2(n+1)\ell_{n+1}^{\alpha}\,. \label{eq2intro}%
\end{equation}

These equations are  transparent in the natural setting where
the group acts on holomorphic functions on the upper half plane $\mathbb{C}%
^{+}=\{z=x+iy\mid y>0\}$. That natural setting is related to $L^{2}%
(\mathbb{R}^{+},t^{\alpha}dt)$ by the Laplace transform $\mathcal{L}_{\nu
}(f)(z)=\int_{\mathbb{R^{+}}}e^{-zt}f(t)~t^{\alpha}dt$ which we derived by a
simple application of the restriction principle (see \cite{doz1,doz2,OOe96}).
The Laplace transform is then used to transfer results from the holomorphic
setting to $L^{2}(\mathbb{R}^{+},t^{\alpha}dt)$. Notice that the automorphism
group of the cone $\mathbb{R}^{+}$ is in this case
\[
H\simeq\left\{  \left(
\begin{array}
[c]{cc}%
e^{t} & 0\\
0 & e^{-t}%
\end{array}
\right)  =\exp\left(  tZ\right)  \mid t\in\mathbb{R}\right\}
\]
where%
\[
Z=\left(
\begin{array}
[c]{cc}%
1 & 0\\
0 & -1
\end{array}
\right)  =\frac{1}{2}(X^{-}-X^{+})~.
\]
The corresponding differential operator is $2tD+\nu$ and we get%
\begin{equation}
\left(  2tD+\nu\right)  \ell_{n}^{\alpha}=(n+1)\ell_{n+1}^{\alpha}%
-(n+\nu-1)\ell_{n-1}^{\alpha} \label{eq4intro}%
\end{equation}
which is, up to normalization of the Laguerre functions, the
recursion relation stated in general in Theorem 7.9 and Remark
7.10 in \cite{doz2}, see equation (\ref{eq3intro}) below.

In our article \cite{doz2} we explore, among other things, a
generalization and extension of these results to the setting where
$\mathbb{R}^{+}$ is replaced by a symmetric cone $H/L=\Omega$,
contained in a Euclidean Jordan algebra $J$, and
$\mathrm{SL}(2,\mathbb{R})$ is replaced by the Hermitian group
$G=\mathrm{Aut}(T(\Omega))_{o}$, where $T(\Omega)=\Omega+iJ$.
Furthermore, the measure $t^{\alpha}dt$ is replaced by
$d\mu_{\nu}(t)=\Delta(t)^{\nu-d/r}dt$, where $\Delta$ is the
determinant function of $J$, $d=\mathrm{dim}(J)$, and $r$ is the
rank of $J$ (see \cite{FK} for the structure theory of Jordan
algebras). In this setting the Laguerre functions form an
orthogonal basis for the space of $L$-invariant functions
$L^{2}(\Omega,d\mu_{\nu})^{L}$. And again, the differential
operators and recurrence relations that they satisfy can be
derived from the action of the complexified Lie algebra
$\mathfrak{g}_{\mathbb{C}}$ of $G$. If
$\mathfrak{g}_{\mathbb{C}}^L$ denotes the $L$-invariant elements
of $\mathfrak{g}_\mathbb{C}$ then $\mathfrak{g}_{\mathbb{C}}^L$
acts
on the space of Laguerre functions.  If $\mathfrak{g}$ is simple then $\dim(\mathfrak{g}%
_{\mathbb{C}}^{L})=3$ and in fact there is a Lie algebra isomorphism%
\[
\varphi:\mathfrak{g}_{\mathbb{C}}^{L}\rightarrow\mathfrak{sl}(2,\mathbb{C}).
\]
We use the same notation as before and let $X^{0}$  be taken from
the center of $\mathfrak{k}$, where $\mathfrak{k}$ is the Lie
algebra of a maximal compact
subgroup $K$ of $G$, containing $L$. If $\sigma:\mathfrak{g}_{\mathbb{C}%
}\rightarrow\mathfrak{g}_{\mathbb{C}}$ is the conjugation with
respect to the real form $\mathfrak{g}$ and $\theta$ is the Cartan
involution such that $G^{\theta}=K$, then we can choose the
isomorphism $\varphi$ such that
$\varphi(\sigma(X))=\overline{\varphi(X)}$ and
$\varphi(\theta(X))=-\varphi (X)^{t}$. We then fix $X^{0}$ such
that
\[
\varphi(X^{0})=\left(
\begin{array}
[c]{cc}%
0 & -i\\
i & 0
\end{array}
\right)  ~.
\]
The operator $\mathrm{ad}(X^{0}):\mathfrak{g}_{\mathbb{C}}%
\rightarrow\mathfrak{g}_{\mathbb{C}}$ has eigenvalues $0,2,$ and $-2$. The
eigenspace corresponding to the eigenvalue $0$ is $\mathfrak{k}_{\mathbb{C}}$
whereas the eigenspace $\mathfrak{p}^{\pm}$ corresponding to the eigenvalue
$\pm2$ is an abelian Lie algebra and $\sigma(\mathfrak{p}^{+})=\mathfrak{p}%
^{-}$.\ In particular, $\mathfrak{p}^{+}\cap\mathfrak{p}^{-}=\{0\}$. We choose
$X^{\mp}\in(\mathfrak{p}^{\pm})^{L}$ such that%
\[
\varphi(X^{-})=\left(
\begin{array}
[c]{cc}%
1 & i\\
i & -1
\end{array}
\right)  ~,\quad\mathrm{and}\quad\varphi(X^{+})=\left(
\begin{array}
[c]{cc}%
-1 & i\\
i & 1
\end{array}
\right)  ~.
\]
Then $Z=\frac{1}{2}(X^{-}-X^{+})$ is in the center of $\mathfrak{h}$. In
\cite{doz2} it was shown, that $Z$ corresponds to the differential operator
\[
E_{\nu}=\nu r+2E,
\]
where $E$ is the \textit{Euler operator} $Ef(x)=\frac{d\,}{dt}f(tx)|_{t=1}$.
Furthermore, the following difference equation was derived:
\begin{equation}
2E_{\nu}\ell_{\mathbf{m}}^{\nu}=\sum_{j=1}^{r}c_{\mathbf{m}}(j)\ell_{\mathbf{m}
+\gamma_{j}}^{\nu}-\sum
_{j=1}^{r}{\binom{\mathbf{m}}{\mathbf{m}-\gamma_{j}}}(\mathbf{m}_{j}-1+\nu-(j-1)
\frac{d}{2})\ell_{\mathbf{m}-\gamma_{j}%
}^{\nu} \label{eq3intro}%
\end{equation}
where $c_{m}(j)$ are certain constants. We notice that this
equation includes the annihilation operator and the creation
operator. It is one of our aims in this paper to separate the
action of the creation and annihilation operators. We therefore
return to the paradigm discussed earlier for
$\mathrm{SL}(2,\mathbb{R})$ and carry out our extension for the
group $\mathrm{SU}(n,n)$ acting on the tube type domain
$T(\Omega)$, where $\Omega$ is the cone of positive definite
$n\times n$ Hermitian matrices and $T(\Omega)$ is $\Omega+ i
\mathrm{Herm}(n)$. Using the Laplace transform we determine the
action of $\mathfrak{g}_\mathbb{C}$ on  $L^{2}(\Omega
,d\mu_{\nu}).$  This calculation is carried out on
$\mathfrak{k}_\mathbb{C}$, where $\mathfrak{k}$ is the Lie algebra
of $K$, on $\mathfrak{p}^+$, and on $\mathfrak{p}^-$.  The result
is:

\medskip

\noindent\textbf{Theorem \ref{th2.1}} \textit{Let} $f\in L^{2}(\Omega,d\mu_{\nu}%
)$\textit{ be a smooth vector.}

\begin{enumerate}
\item Let $%
\begin{pmatrix}
a & b\\
b & a
\end{pmatrix}
\in\mathfrak{k}_{\mathbb{C}}$. Then
\[
\lambda_{\nu}
\begin{pmatrix}
a & b\\
b & a
\end{pmatrix}
f(s)=\mathrm{tr}(-s\nabla b\nabla+(sa-as-\nu b)\nabla + b s)f(s)
\]
\item Let $%
\begin{pmatrix}
x & x\\
-x & -x
\end{pmatrix}
\in\mathfrak{{p}^{+}}$. Then
\[
\lambda_{\nu}
\begin{pmatrix}
x & x\\
-x & -x
\end{pmatrix}
f(s)=\mathrm{tr}(s\nabla x\nabla+(\nu x+sx+xs)\nabla+(\nu
x+sx))f(s).
\]
\item Let $%
\begin{pmatrix}
x & -x\\
x & -x
\end{pmatrix}
\in\mathfrak{{p}^{-}}$. Then
\[
\lambda_{\nu}%
\begin{pmatrix}
x & -x\\
x & -x
\end{pmatrix}
f(s)=\mathrm{tr}(-s\nabla x\nabla+(-\nu x+sx+xs)\nabla+(\nu
x-sx))f(s).
\]

\end{enumerate}

Specializing to $\mathfrak{g}_\mathbb{C}^L=\operatorname{span}
(X_\circ,X^-,X^+)$ gives the analogue to Laguerres differential
operator, the annihilation operator, and the creation operator.
Specifically, we get

\medskip

\noindent\textbf{Theorem \ref{thm5.1}} \textit{The Laguerre
functions are related by the following recursion relations:}

\begin{enumerate}

\item $\mathrm{tr}(-s\nabla\nabla-\nu \nabla+s)\ell_{\mathbf{m}}^{\nu}=(r\nu+|\mathbf{m}|)\ell_{\mathbf{m}}^{\nu}.$

\item $\frac{1}{2}\mathrm{tr}(s\nabla\nabla+(\nu I+2s)\nabla+(\nu I+s)
)\ell_{\mathbf{m}}^{\nu}(s)=-\sum_{j=1}^{r}%
\begin{pmatrix}
\mathbf{m}\\
\mathbf{m-\gamma_{j}}%
\end{pmatrix}
(m_{j}-1+\nu-(j-1)\frac{d}{2})\ell_{\mathbf{m-\gamma_{j}}}^{\nu} $

\item $\frac{1}{2}\mathrm{tr}(-s\nabla\nabla+(-\nu I+2s)\nabla+(\nu I-s)
)\ell_{\mathbf{m}}^{\nu}=\sum_{j=1}^{r}c_{\mathbf{m}}(j)\ell_{\mathbf{m+\gamma
_{j}}}^{\nu}.$
\end{enumerate}
\medskip

Several aspects of these formulas suggest an interpretation
involving the structure theory of Jordan algebras. We have not
pursued that here.

Our realization of the highest weight representations
$\lambda_{\nu}$ in $L^{2}(\Omega,d\mu_{\nu})$ has an interesting
connection to the set of Whittaker vectors of the representation
$\lambda_{\nu}$, see \cite{K00,W01}. Recall that
$\mathfrak{z}(\mathfrak{h})=\mathbb{R}Z.$ Let
\begin{equation}
\mathfrak{n}=\{X\in\mathfrak{g}\mid\lbrack Z,X]=2X\}\,.
\end{equation}
Then $\mathfrak{n}$ is an abelian Lie algebra, and $\mathfrak{h}%
\oplus\mathfrak{n}=:\mathfrak{p}_{\mathrm{max}}$ is a maximal parabolic
subalgebra in $\mathfrak{g}$. Let $P=\tilde{H}N$ be the corresponding maximal
parabolic subgroup and notice that $\tilde{H}_{o}=H$. Let us recall here the
construction of the isomorphism $\mathcal{H}_{\nu}(\Omega)\simeq L^{2}%
(\Omega,d\mu_{\nu})$ where $\mathcal{H}_{\nu}(\Omega)\subset\mathcal{O}%
(T(\Omega))$ is the usual Hilbert space realization of the highest
weight representation $\pi_{\nu}$. Then the restriction map
$R(f):=f|_{\Omega}$ is injective and intertwines in an obvious way
the $\tilde{H}$-action on both sides. Now it is a well known fact
that the cone $\Omega$ has a natural realization in the vector
space $\mathfrak{n}\simeq N$ and therefore this map also
intertwines in a natural way the $N$-action on both sides. It then
follows that $(\pi_{\nu},\mathcal{H}_{\nu})$ can be realized as
$(\tilde{\lambda}_{\nu},L^{2}(\Omega,dt))$ where we think now of
$L^{2}(\Omega,dt)$ as a subspace of $L^{2}(N)$. The isomorphism
$L^{2}(\Omega,d\mu_{\nu})\simeq L^{2}(\Omega,dt)$ is here given by
\begin{equation}
f\mapsto\Delta^{(\nu-d/r)/2}f
\end{equation}
In our example of $G=\mathrm{SL}(2,\mathbb{R})$ we have $Z=\left(
\begin{array}
[c]{cc}%
1 & 0\\
0 & -1
\end{array}
\right)  $,
\[
N=\left\{  \left(
\begin{array}
[c]{cc}%
1 & x\\
0 & 1
\end{array}
\right)  \mid x\in\mathbb{R}\right\}  \simeq\mathbb{R}\,.
\]
and
\[
P=\pm\left\{  \left(
\begin{array}
[c]{cc}%
a & x\\
0 & a^{-1}%
\end{array}
\right)  \mid a>0,x\in\mathbb{R}\right\}  \,.
\]
Furthermore, the isomorphism
$L^{2}(\mathbb{R}^{+},t^{\alpha}dt)\simeq
L^{2}(\mathbb{R}^{+},dt)$ is simply
\[
f\mapsto(t\mapsto t^{\alpha/2}f(t))\,.
\]
Let%
\[X_N=\left(
\begin{array}
[c]{cc}%
0 & 1\\
0 & 0
\end{array}
\right)  =\frac{-i}{2}\left(  \frac{1}{2}(X^{+}+X^{-})-X^{0}\right)  ~.
\]
The set of equations (\ref{slaction1}-\ref{slaction3})\ gives us
then the
following action of the Lie algebra on functions on $L^{2}(\mathbb{R}^{+}%
,dt)$:%
\begin{align*}
Z\cdot f(t)  &  =2tf^{\prime}(t)+f(t)\\
X_{N}\cdot f(t)  &  =-tif(t)\\
(iX_{0})\cdot f(t)  &  =i\left(  tf^{\prime\prime}(t)+f^{\prime}(t)-\left(
\frac{\alpha^{2}}{4t}+t\right)  f(t)\right)
\end{align*}
These are exactly the equations in \cite{K00}. The motivation in
\cite{K00} is  different though. There, as well as in \cite{W01},
the differential equations satisfied by the Laguerre functions
were used to reconstruct the representation $\lambda_{\nu}$
whereas we have used the representation $\lambda_{\nu}$ to derive
the differential equations that the Laguerre functions satisfy.
It
is now an interesting task to generalize this for the universal covering of
the group $\mathrm{SU}(n,n)$. We remark the following. First notice that $r/d=1/n$.
The isomorphism
$L^{2}(\Omega,\Delta(t)^{\nu-1/n}\, dt)\simeq L^{2}(\Omega,dt)$ is again given by%
\[
f\mapsto\Delta^{\frac{n\nu-1}{2n}}f\, .
\]
Now one should use Theorem \ref{th2.1} and Proposition \ref{p3.2}, which gives
us the derivatives of $\Delta^{\frac{n\nu-1}{2n}}$, to derive
the the action of $\mathfrak{su}(n,n)$ on $L^{2}(\Omega,dt)$.

After we submitted this article the paper of Ricci and Vignati
\cite{RV} came to our attention.  There a system of differential
operators are defined that are diagonal on the Laguerre basis with
different eigenvalues.  This system generalizes the classical
Laguerre operator, however, the raising and lowering operators are
not discussed.

\section{$\mathrm{SU}(n,n)$ and its Lie Algebra}
\label{sunn}

 We keep the
notation from the introduction but specialize to the case $G=\
\mathrm{SU}(n,n)$. The standard definition of $\mathrm{SU}(n,n)$
as found on page 444 of \cite{helg} acts naturally on the
generalized unit disk. Our definition below is suited for the
right half plane action and is conjugate to the standard version
by the Cayley transform. The group $\mathrm{SU}(n,n)$ is thus
defined as follows: Let $J=\begin{pmatrix}
0 & 1\\
1 & 0
\end{pmatrix}
$. Then
\[
\mathrm{SU}(n,n)=\left\{  g\in\mathrm{SL}(2n,\mathbb{C})\mid gJg^{\ast}=J \right\}  .
\]
We frequently write $g\in SU(n,n)$ in block form as $g=%
\begin{pmatrix}
A & B\\
C & D
\end{pmatrix}\, .$ We then have the following relations on the entries:
\begin{align*}
AB^{\ast}  +BA^{\ast}&=0  \qquad B^*D+D^*B=0\\
CD^{\ast}  +DC^{\ast}&=0 \qquad A^*C+C^*A=0\\
AD^{\ast}+BC^{\ast}  &  =1\qquad A^*D+C^*B=1\\
DA^{\ast}+CB^{\ast}  &  =1\qquad B^*C+D^*A=1.\\
\end{align*}

The cone $\Omega$ in $J=\mathrm{Herm}(n)$ can be identified as the
set of positive definite Hermitian matrices. The group
$\mathrm{SU}(n,n)$ acts on $T(\Omega)=\Omega+iJ$ by linear
fractional transformations:
\[
\left(
\begin{array}
[c]{cc}%
A & B\\
C & D
\end{array}
\right)  \cdot T=(AT+B)(CT+D)^{-1}.
\]

Let $\Gamma_{\Omega}(\nu )$ be the \textit{Gindikin-Koecher Gamma
function} associated with the cone $\Omega$, see \cite{FK,SG}:
\[\Gamma_{\Omega}(\nu) =\int_{\Omega}e^{-\mathrm{tr}(x)}\det(x)^{\nu -n}dx\, .\]
Assume that $\nu >2n-1$. Let $\mathcal{H}_{\nu}(T(\Omega))$ be the space of holomorphic functions
$F:T(\Omega)\to \mathbb{C}$ such that
\begin{eqnarray}\label{eq-Hnu}
\|F\|^2:= \alpha_\nu\int_{T(\Omega)}|F(x+iy)|^2\det(y)^{\nu-2n}\,
dxdy <\infty
\end{eqnarray}
where
\[\alpha_\nu=\frac{2^{n\nu}}{(4\pi)^{n^2}\Gamma_{\Omega}(\nu -n)}\, .\]
Then $\mathcal{H}_{\nu}(T(\Omega))$ is a Hilbert space and  the
group $\mathrm{SU}(n,n)$ acts unitarily and irreducibly on
$\mathcal{H}_{\nu}(T(\Omega))$ by
\[
\pi_{\nu}(g)F(z)=J(g^{-1},z)^{\frac{\nu}{2n}}f(g^{-1}\cdot z),
\]
where $J(%
\begin{pmatrix}
A & B\\
C & D
\end{pmatrix} ,z)=\det(Cz+D)^{-2n}$, (see  section 2.8 of
\cite{doz2}). Let us also mention, that the Hilbert space
structure and the unitary action extends analytically to the
domain $\nu > n-1$, \cite{RV76,W79}.

The following table summarizes most aspects of the specialization
from the introduction to $\mathrm{SU}(n,n)$, some of which we have
already used:
\renewcommand*{\arraystretch}{2.0}

\begin{center}%
\begin{tabular}
[c]{|c|c|}\hline $G$ & $\mathrm{SU}(n,n)$\\\hline\hline $a$ & 2
\\\hline $d$ & $n^{2}$\\\hline $r$ & $n$\\\hline $p$ & $2n$\\\hline
$K$ & \renewcommand*{\arraystretch}{1.0} $%
\begin{pmatrix}
A & B\\
B & A\
\end{pmatrix}
\in\mathrm{SU}(n,n):\quad A\pm B\in\mathrm{U}(n)$
\renewcommand*{\arraystretch}{2.0}\\\hline
$H$ & \renewcommand*{\arraystretch}{1.0}$%
\begin{pmatrix}
A & 0\\
0 & (A^{\ast})^{-1}%
\end{pmatrix}
\in\mathrm{SL}(2n,\mathbb{C}):\quad A\in\mathrm{GL}(n,\mathbb{C}%
)$\renewcommand*{\arraystretch}{2.0}\\\hline
$L$ & \renewcommand*{\arraystretch}{1.0} $%
\begin{pmatrix}
A & 0\\
0 & A\
\end{pmatrix}
\in\mathrm{SL}(2n,\mathbb{C}):\quad A\in\mathrm{SU}(n)$
\renewcommand*{\arraystretch}{2.0}\\\hline
$J$ & Herm(n): the $n\times n$ Hermitian matrices\\\hline
$\Omega$ & positive definite matrices in $J$\\\hline
$\Delta$ & $\det$\\\hline
\end{tabular}

\end{center}

\renewcommand*{\arraystretch}{1.0}

\subsection{The Lie Algebra}

The Lie algebra for\textrm{\ }$\mathrm{SU}(n,n)$ is
\[
\mathfrak{g}=\mathfrak{su}(n,n)=\left\{
\begin{pmatrix}
a & b\\
c & d
\end{pmatrix}\in \mathfrak{sl}(2n,\mathbb {C})
\mid a=-d^{\ast},b=-b^{\ast},c=-c^{\ast}\right\} .
\]
Its complexification, $\mathfrak{g}_{\mathbb{C}}$, is the set
$\mathfrak{sl}(2n,\mathbb{C})$. Important Lie subalgebras are
$\mathfrak{k}$, the Lie algebra of $K$, its complexification
$\mathfrak{k}_{\mathbb{C}}$, and $\mathfrak{p}^{+}$ and
$\mathfrak{p}^{-}$, which act by annihilation  and creation
operators in $\mathcal{H}_{\nu}(T(\Omega))$. They are given as
follows:
\[
\mathfrak{k}=\left\{
\begin{pmatrix}
a & b\\
b & a
\end{pmatrix}
\mid a=-a^{\ast},\;b=-b^{\ast},\text{ and
}\operatorname{tr}(a)=0\right\}  ,
\]%
\[
\mathfrak{k}_{\mathbb{C}}=\left\{
\begin{pmatrix}
a & b\\
b & a
\end{pmatrix}
\mid\operatorname{tr}(a)=0\right\}  ,
\]%
\[
\mathfrak{p}^{+}=\left\{
\begin{pmatrix}
x & x\\
-x & -x
\end{pmatrix}
\mid x=x^{\ast}\right\}  ,
\]
and
\[
\mathfrak{p}^{-}=\left\{
\begin{pmatrix}
x & -x\\
x & -x
\end{pmatrix}
\mid x=x^{\ast}\right\}  .
\]

\section{The actions of the subalgebras $\mathfrak{k}, \; \mathfrak{p}%
^{-},\text{ and }\mathfrak{p}^{+}$}

\label{sunnaction} In this section we give the action of the complexification,
$\mathfrak{g}_{\mathbb{C}}$, of the Lie algebra of $G=\mathrm{SU}(n,n)$ on
$\mathcal{H}_{\nu}(T(\Omega))$. The proof is a straightforward calculation and
will not be given here. Recall though the basic definitions. Let
$X\in\mathfrak{g}$ and let $F\in\mathcal{H}_{\nu}(T(\Omega))$ be a smooth
vector. By definition
\begin{align*}
\pi_{\nu}(X)F(z) & =\frac{d\, }{dt}\pi_{\nu}(\exp tX)F(z)|_{t=0}\\
& =\frac{d\, }{dt} J(\exp-tX,z)^{\frac{\nu}{2n}}F(\exp(-tX)\cdot
z)|_{t=0}.
\end{align*}
The representation of $\mathfrak{g}_{\mathbb{C}}$ is defined by $\pi_{\nu
}(X+iY)=\pi_{\nu}(X)+i\pi_{\nu}(Y)$. But since functions in $\mathcal{H}_{\nu
}(T(\Omega))$ are holomorphic we get that $\pi_{\nu}(X)$ is given by the same
formula for $X\in\mathfrak{g}_{C}$.

\begin{thm}
\label{aonH} For $w\in V$ let $\delta(w)$ denote the
non-normalized directional derivative in the direction $w$. Let
$F$ be a $C^{\infty}$-vector in $\mathcal{H}_{\nu}(T(\Omega))$.

\begin{enumerate}
\item Let $X=%
\begin{pmatrix}
a & b\\
b & a
\end{pmatrix}
\in\mathfrak{k}_{\mathbb{C}}.$ Then
\[
\pi_{\nu}(X)F(z)=\nu\operatorname{tr}(bz)F(z)+\delta(zbz+za-az-b)F(z).
\]

\item Let $X=%
\begin{pmatrix}
x & x\\
-x & -x
\end{pmatrix}
\in\mathfrak{p}^{+}. $ Then
\[
\pi_{\nu}(X)F(z)=-\nu\operatorname{tr}({x(1+z)})F(z)-\delta( {(1+z)x(1+z)}%
)F(z).
\]

\item Let $X=%
\begin{pmatrix}
x & -x\\
x & -x
\end{pmatrix}
\in\mathfrak{p}^{-}. $ Then
\[
\pi_{\nu}(X)F(z)=\nu\operatorname{tr}({x(z-1)})F(z)+\delta( {(z-1)x(z-1)}%
)F(z).
\]

\end{enumerate}
\end{thm}

\section{The gradient}

\label{gradient} Our main results are expressed in terms of the gradient of
complex valued function. This section establishes notation and well known results.

Let $V$ be a finite dimensional complex vector space with inner product
$( \cdot\mid \cdot)$. Suppose $J$ is a real form on which the restriction of the
inner product is a real inner product. Let $f:J\rightarrow\mathbb{R}$ be a
differentiable function, i.e., all (nonnormalized) directional derivatives
$D_{u},u\in J$ , exist. For $s\in J$, we define $\nabla f(s)\in J$ by the
formula
\[
(\nabla f(s)\mid u)=D_{u}f(s).
\]
For a $\mathbb{C}$-valued function $f=f_{1}+if_{2}$ we define $\nabla f=\nabla
f_{1}+i\nabla f_{2}.$ For $z=x+iy\in V$, $x,y\in J$, we define $D_{z}%
=D_{x}+iD_{y}$. The following proposition is for the most part a
consequence of these definitions.

\begin{prop} Let the notation be as above. Then the following holds:
\begin{enumerate}
\item Let $f:J\rightarrow\mathbb{C}$ be differentiable. Then
\[
(\nabla f(s)\vert z)=D_{\bar{z}} f(s),
\]
where $\bar{z}$ is conjugation with respect to $J$.

\item If $f:J\to \mathbb{C}$ is the restriction of a  $\mathbb{C}$-linear map $F:V\to \mathbb{C}$ to $J$, then
$(\nabla f(s)\mid w)=F(\bar{w})$. In particular if $z,w\in V$ and $f(s)
=(z\mid s)=(s\mid \bar{z})$ then $(\nabla f(s)\mid w)
=([\nabla(z\mid s)]\mid w)=(z\mid w)$.

\item If $f$ is the restriction of a holomorphic function on $V$ to $J$
then $D_z=\delta_z$. I.e.
\[
D_{z}f(s)=\frac{d}{dt}f(s+tz)|_{t=0},
\]
is the (non-normalized) directional derivative in $V$.

\item For $z,w\in V$ $\ $we have $D_{w}(z\mid s)=(z\mid\bar{w}).$
\end{enumerate}
\end{prop}

Suppose $z=x+iy$, $x,y\in J,$ $x\ne 0$, $y\ne0$. Then $z$ and
$\bar{z}$ are independent over $\mathbb{C}$, and
$\mathrm{span}_{\mathbb{C}}\{z,\bar z\} =
\mathrm{span}_{\mathbb{C}}\{x,y\}.$ On the other hand, if $z=x$
then $\mathrm{span}_{\mathbb{C}}\{z\}=
\mathrm{span}_{\mathbb{C}}\{x\}$ and if
$z=iy$ then $\mathrm{span}_{\mathbb{C}}\{z\}= \mathrm{span}_{\mathbb{C}%
}\{y\}.$ It follows from these observations that there is an orthonormal basis
of $V$ of the form
\[
\{e_{1},\ldots, e_{n}, z_{1}, \bar{z_{1}}, \ldots, z_{k}, \bar{z_{k}}\},
\]
where $e_{1},\ldots, e_{n} \in J$ and
\[
\{e_{1},\ldots,e_{n},x_{1},\ldots,x_{k},y_{1},\ldots y_{k}\},
\]
where $z_{j}=x_{j} + iy_{j}, \;j=1,\cdots,k$, forms a basis of J. We call such
a basis of $V$ \textbf{compatible} with $J$.

Now let $J$ be the $\mathbb{R}$-vector space of $n\times n$ Hermitian
symmetric spaces. I.e. $J=\{s\in M_{n,n}(\mathbb{C}):s=s^{\ast}\}$. On
$V=M_{n,n}(\mathbb{C})$ there is a complex inner product given by
$(z_{1}|z_{2})=\mathrm{tr}(z_{1}z_{2}^{\ast})$ which, when restricted to $J$,
is a real scalar product. Notice also that conjugation with respect to the
real form J is given by the adjoint ${}^\ast$, i.e., $\bar{v}=v^{\ast}$. It is
convenient to pick a specific basis for the calculations in the next section.
Let $E_{i,j}$ be the $n\times n$ matrix with $1$ in the $(i,j)$-entry and
$0$'s elsewhere. Then the collection $\{E_{i,j}\}$ is a basis of $V$
compatible with $J$. For, observe that $E_{i,i}^{\ast}=E_{i,i}$ and
$E_{i,j}^{\ast}=E_{j,i}$, $i\neq j$. Let $I_{i}=E_{i,i},\;J_{i,j}%
=E_{i,j}+E_{j,i},\;\text{ and }K_{i,j}=iE_{i,j}-iE_{j,i}.$ (There
should be no confusion between the
index $i$  a positive integer and the complex number $i=\sqrt{-1}$.) Then $E_{i,j}%
=\frac{1}{2}(J_{i,j}-iK_{i,j}).$  Furthermore
\[
\mathcal{B}=\{I_{i},\;\frac{1}{2}J_{i,j},\;\frac{1}{2}K_{i,j}\},
\]
where $i<j$, is a corresponding basis of $J$. Let $D_{i,j}%
=D_{E_{i,j}}=D_{\frac{1}{2}(J_{i,j}-iK_{i,j})}.$ Using this notation it
follows that
\[
(\nabla)_{i,j}f:=(\nabla f)_{i,j}=(\nabla f|E_{i,j})=D_{\bar{E}_{i,j}%
}f=D_{j,i}f.
\]
Furthermore, for $w\in V$,
\[
D_{w}f=\sum w_{i,j}D_{i,j}f=\sum w_{i,j}(\nabla)_{j,i}f=\mathrm{tr}%
(w\nabla)f,
\]
where $\mathrm{tr}$ is the usual trace functional.

\begin{prop}\label{p3.2}
\begin{enumerate}
\item[]
\item For $1\leq i,j,k,l\leq n$ we have
\[
D_{i,j}s_{k,l}=\delta_{i,k}\delta_{j,l}.
\]

\item Suppose $w\in V$ and $s\in J$ is invertible. Then
\[
D_{w}\det(s)^{m}=m\det(s)^{m}\mathrm{tr}(s^{-1}w)\, .
\]

\item For $1\leq i,j\leq n$ we have
\[
D_{i,j}\det(s)^{m}=m\det(s)^{m-1}\operatorname{cof}_{i,j}(s),
\]
where $\operatorname{cof}_{i,j}(s)$ denotes the $(i,j)$-cofactor of $s$
\end{enumerate}
\end{prop}

\begin{proof}
(1) is clear and (3) follows from (2). Since $\det$ is the restriction of a
holomorphic function we have
\[
D_{w}\det(s)=\frac{d}{dt}\det(s)\det(1+ts^{-1}w)|_{t=0}=\det(s)\mathrm{tr}%
(s^{-1}w).
\]
Thus
\[
D_{w}\det(s)^{m}=m\det(s)^{m-1}D_{w}\det(s)=m\det(s)^{m}\mathrm{tr}(s^{-1}w).
\]

\end{proof}

\begin{prop}
Suppose that $f,g\in L^2(\Omega, d\mu_{\nu})$ are smooth and $f$
vanishes on the boundary of the cone $\Omega$. Let $1\le i,j \le
n$. Then
\begin{enumerate}
\item
$$ \int_{\Omega}D_{i,j} f(s) g(s) ds=-\int_{\Omega}f(s) D_{i,j}
g(s) ds. $$
\item
$$ \int_{\Omega} e^{-(z\vert s)}z_{i,j} f(s) ds = \int_{\Omega} e^{-(z\vert
s)}D_{j,i} f(s) ds. $$
\end{enumerate}
\end{prop}

\begin{proof}
(1) is Stokes Theorem and (2) follows from (1) and the fact that
$D_{j,i}e^{-(z\vert s)}=-e^{-(z\vert s)}z_{i,j}.$
\end{proof}

\section{The action of $\mathfrak{g}_{\mathbb{C}}$ on $L^{2}(\Omega,d\mu_{\nu
})$}

In section \ref{sunnaction} we determined the algebraic action of
$\mathfrak{g}_{\mathbb{C}}$ on $\mathcal{H}_{\nu}(T(\Omega))$. Recall the
measure $d\mu_{\nu}(t)=\det(t)^{\nu-n}dt~$and the inner product $(x\mid
y)=\mathrm{tr}(xy)$. The Laplace transform
\[
\mathcal{L}_{\nu}:L^{2}(\Omega,d\mu_{\nu})\rightarrow\mathcal{H}_{\nu
}(T(\Omega)),
\]
given by the formula
\[
\mathcal{L}_{\nu}(f)(z)=\int_{\Omega}e^{-(z|s)}f(s)~\det(s)^{\nu-n}%
ds=\int_{\Omega}e^{-(z\mid s)}~f(s)~d\mu_{\nu}(s),
\]
defines an isomorphism between  $L^{2}(\Omega,d\mu_{\nu})$ and
$\mathcal{H}_{\nu}(T(\Omega))$. We can thus define an equivalent
action,  denoted $\lambda_{\nu}$, of $\mathfrak{g}_{\mathbb{C}}$
on $L^{2}(\Omega,d\mu_{\nu})$. This action is given in the
following theorem.

\begin{thm}
\label{th2.1}Let $f\in L^{2}(\Omega,d\mu_{\nu})$ be a smooth vector.
\begin{enumerate}
\item Let $%
\begin{pmatrix}
a & b\\
b & a
\end{pmatrix}
\in\mathfrak{k}_{\mathbb{C}}$. Then
\[
\lambda_{\nu}
\begin{pmatrix}
a & b\\
b & a
\end{pmatrix}
f(s)=\mathrm{tr}(-s\nabla b\nabla+(s a-as-\nu b)\nabla +b s)f(s)
\]
\item Let $%
\begin{pmatrix}
x & x\\
-x & -x
\end{pmatrix}
\in\mathfrak{{p}^{+}}$. Then
\[
\lambda_{\nu}
\begin{pmatrix}
x & x\\
-x & -x
\end{pmatrix}
f(s)=\mathrm{tr}(s\nabla x\nabla+(\nu x+sx+xs)\nabla+(\nu
x+sx))f(s).
\]
\item Let $%
\begin{pmatrix}
x & -x\\
x & -x
\end{pmatrix}
\in\mathfrak{{p}^{-}}$. Then
\[
\lambda_{\nu}%
\begin{pmatrix}
x & -x\\
x & -x
\end{pmatrix}
f(s)=\mathrm{tr}(-s\nabla x\nabla+(-\nu x+sx+xs)\nabla+(\nu
x-sx))f(s).
\]

\end{enumerate}
\end{thm}

The proof of this theorem is somewhat long and certainly tedious. We will
provide the details for the subalgebra $\mathfrak{p}^{+}$. The actions for the
other algebras are done in a similar manner. Throughout we will use without
reference the results of section \ref{gradient}.

\subsection{(2) The proof for the $\mathfrak{p}^{+}$- action}

Let $X=%
\begin{pmatrix}
x & x\\
-x & -x
\end{pmatrix}
\in\mathfrak{p}^{+}$ and let $f\in L^{2}(\Omega,d\mu_{\nu})$ be a
smooth vector. We use the Laplace transform $\mathcal{L}_{\nu}$ to
transfer the action of $\mathfrak{p}^{+}$ on
$L^{2}(\Omega,d\mu_{\nu})$. Thus we seek a formula for
$\lambda_{\nu}(X)$ such that
\[
\pi_{\nu}(X)\mathcal{L}_{\nu}f(z)=\mathcal{L}_{\nu}(\lambda_{\nu}(X)f)(z).
\]
Recall from Theorem \ref{aonH} that for $F$ a smooth vector in
$\mathcal{H}_{\nu}(T(\Omega))$ we have
\[
\pi_{\nu}(X)F(z)=-\nu\operatorname{tr}( {x(1+z)})F(z)-\delta(
{(1+z)x(1+z)}) F(z).
\]
Our first calculation will begin with
$-\delta((z+1)x(z+1))\mathcal{L}_{\nu }f(z)$. Let
$m=\nu-d/r=\nu-n.$

\begin{lem}
\
\begin{align}
-\delta((z+1)x(z+1))\mathcal{L}_{\nu}f(z)  &  = -\delta((z+1)x(z+1))\int
e^{-(z|s)}f(s)\det(s)^{m}~ds\nonumber\\
&  =\int e^{-(z|s)}((z+1)x(z+1)|s)f(s)~\det(s)^{m}~ds\nonumber\\
&  =\int e^{-(z|s)}(x|s)f(s)~\det(s)^{m}~ds\tag{A}\\
&  \quad +\int e^{-(z|s)}(xz|s)f(s)~\det(s)^{m}~ds\tag{B}\\
&  \quad +\int e^{-(z|s)}(zx|s)f(s)~\det(s)^{m}~ds\tag{C}\\
&  \quad +\int e^{-(z|s)}(zxz|s)f(s)~\det(s)^{m}~ds. \tag{D}%
\end{align}
\end{lem}

\begin{lem}
[B]\label{B}
\[
\int
e^{-(z|s)}(xz|s)f(s)\det(s)^{m}ds=(n+m)\mathrm{tr}(x)\mathcal{L}_{\nu
}f(z) +\int e^{-(z|s)}\mathrm{tr}(sx\nabla)f(s)\det(s)^{m}\, ds.
\]

\end{lem}

\begin{proof}
First observe that $(xz|s)=\sum_{i,j}(xz)_{i,j}s_{i,j}^{\ast}=\sum
_{i,j,k}x_{i,k}z_{k,j}s_{j,i}.$
\begin{eqnarray*}
B  &  =&\sum_{i,k,j}\int e^{-(z|s)}x_{i,k}D_{j,k}\left(
s_{j,i}f(s)\det
(s)^{m}\right)  ~ds\\
&  =&\sum_{i,k,j}\int e^{-(z|s)}x_{i,k}\delta_{k,i}f(s)\det(s)^{m}~ds\\
&  &+\sum_{i,k,j}\int e^{-(z|s)}x_{i,k}s_{j,i}(D_{j,k}f)(s)\det(s)^{m}~ds\\
&  &+\sum_{i,k,j}\int e^{-(z|s)}x_{i,k}s_{j,i}f(s)m\det(s)^{m-1}%
\operatorname{cof}_{j,k}(s)~ds\\
&  =&(n+m)\mathrm{tr}(x)\mathcal{L}_{\nu}f(z) +\sum_{j,k}\int e^{-(z|s)}%
(sx)_{j,k}(D_{j,k}f)(s)\det(s)^{m}~ds\\
&  =&(n+m)\mathrm{tr}(x)\mathcal{L}_{\nu}f(z) +\int e^{-(z|s)}\mathrm{tr}%
(sx\nabla)f(s)\det(s)^{m}~ds.
\end{eqnarray*}
\end{proof}

\begin{lem}
[C]\label{C}
\[
\int
e^{-(z|s)}(zx|s)f(s)\det(s)^{m}ds=(n+m)\mathrm{tr}(x)\mathcal{L}_{\nu
}f(z) +\int e^{-(z|s)}\mathrm{tr}(xs\nabla)f(s)\det(s)^{m}~ds.
\]
\end{lem}
\begin{proof} This proof is very similar to the proof of lemma
\eqref{B}. It is omitted.
\end{proof}

Putting B and C together gives
\begin{lem}
\[
B+C=2(n+m)\mathrm{tr}(x)\mathcal{L}_{\nu}f (z)
+\mathcal{L}_{\nu}(tr((xs+sx)\nabla )f)(z).
\]
\end{lem}

We now calculate D and obtain
\begin{lem}
[D]\label{D}
\begin{eqnarray*}
D& =& (n+m)\mathrm{tr}(zx)\mathcal{L}_{\nu}f (z)\\
& &+(n+m)\mathcal{L}_{\nu
}(\mathrm{tr}(x\nabla)f)(z)
\\
&& +\mathcal{L}_{\nu}(\mathrm{tr}(s\nabla
x\nabla)f)(z)\, .
\end{eqnarray*}
\end{lem}

\begin{proof}
First observe that
\[
(zxz|s)=\sum_{i,j,k,l}z_{i,k}x_{k,l}z_{l,j}s_{j,i}.
\]
Thus
\begin{eqnarray*}
D  &  =&\sum_{i,k,j,l}\int e^{-(z|s)}z_{i,k}z_{l,j}s_{j,i}x_{k,l}%
f(s)\det(s)^{m}~ds\\
&  =&\sum_{i,k,j,l}x_{k,l}z_{i,k}\int
e^{-(z|s)}D_{j,l}(s_{j,i}f(s)\det
(s)^{m})~ds\\
&  =&\sum_{i,k,j,l}x_{k,l}z_{i,k}\int
e^{-(z|s)}\delta_{l,i}f(s)\det
(s)^{m}~ds\\
&  &+\sum_{i,k,j,l}x_{k,l}z_{i,k}\int
e^{-(z|s)}s_{j,i}(D_{j,l}f)(s)\det
(s)^{m}~ds\\
&  &+\sum_{i,k,j,l}x_{k,l}z_{i,k}\int e^{-(z|s)}s_{j,i}f(s)m\det(s)^{m-1}%
\operatorname{cof}_{j,l}(s)~ds\\
&  =&(n+m)\sum_{k,l}x_{k,l}z_{l,k}\int e^{-(z|s)}f(s)\det(s)^{m}~ds\\
&  &+\sum_{i,k,j,l}x_{k,l}z_{i,k}\int
e^{-(z|s)}s_{j,i}(D_{j,l}f)(s)\det
(s)^{m}~ds\\
&  =&(n+m)\mathrm{tr}(zx)\mathcal{L}_{\nu}f(z)+\sum_{i,k,j,l}x_{k,l}\int
e^{-(z|s)}D_{k,i}(s_{j,i}(D_{j,l}f)(s)\det(s)^{m})~ds\\
&  =&(n+m)\mathrm{tr}(zx)\mathcal{L}_{\nu}f(z) \\
&  &+\sum_{i,k,j,l}x_{k,l}\int
e^{-(z|s)}\delta_{k,j}(D_{j,l}f)(s)\det
(s)^{m}~ds\\
&  &+\sum_{i,k,j,l}x_{k,l}\int
e^{-(z|s)}s_{j,i}(D_{k,i}D_{j,l}f)(s)\;\det
(s)^{m}~ds\\
&  &+\sum_{i,k,j,l}x_{k,l}\int
e^{-(z|s)}s_{j,i}(D_{j,l}f)(s)\;m\det
(s)^{m-1}\operatorname{cof}_{k,i}(s)~ds\\
&  =&(n+m)\mathrm{tr}(zx)\mathcal{L}_{\nu}f(z) \\
&  &+(n+m)\int e^{-(z|s)}\sum_{k,l}x_{k,l}(D_{k,l}f)(s)\det(s)^{m}~ds\\
&  &+\int
e^{-(z|s)}\sum_{i,k,j,l}s_{j,i}(D_{k,i}x_{k,l}(D_{j,l}f))(s)\det
(s)^{m}~ds\\
&  =&(n+m)\mathrm{tr}(zx)\mathcal{L}_{\nu}f(z) \\
&  &+(n+m)\mathcal{L}_{\nu}(\mathrm{tr}(x\nabla)f)(z) \\
&  &+\mathcal{L}_{\nu}(\mathrm{tr}(s\nabla x\nabla)f)(z)\, .
\end{eqnarray*}
\end{proof}

By putting together lemmas (4.2)-(4.6) we obtain
\begin{lem}%
\begin{eqnarray*}
-\delta((z+1)x(z+1))\mathcal{L}_{\nu}f(z)
&  =&\mathcal{L}_{\nu} (\mathrm{tr}(sx)f)(z)+2(n+m)\mathrm{tr}(x)\mathcal{L}_{\nu}f(z) \\
& & +\mathcal{L}_{\nu}(D_{xs+sx}f)(z)+(n+m)\mathrm{tr}(zx)\mathcal{L}_{\nu}(f)(z) \\
&
&+(n+m)\mathcal{L}_{\nu}(\mathrm{tr}(x\nabla)f)(z)+\mathcal{L}_{\nu
}(\mathrm{tr}(s\nabla x\nabla)f)(z).
\end{eqnarray*}
\end{lem}

We thus obtain
\begin{eqnarray*}
\pi_{\nu}(%
\begin{pmatrix}
x & x\\
-x & -x
\end{pmatrix}
\mathcal{L}_{\nu}(f)(z)
&  =&-\nu \mathrm{tr}(x(z+1))\mathcal{L}_{\nu}f(z)-\delta
((z+1)x(z+1)\mathcal{L}_{\nu}f(z)\\
&  =&-\nu\mathrm{tr}(xz)\mathcal{L}_{\nu}f (z)-\nu\mathrm{tr}(x)\mathcal{L}_{\nu
}f(z) \\
&&+\mathcal{L}_{\nu}(\mathrm{tr}(sx)f)(z) +2\nu\mathrm{tr}(x)\mathcal{L}_{\nu}f(z)\\
&  &+\mathcal{L}_{\nu}(\mathrm{tr}((sx+xs)\nabla)f)(z) + \nu\mathrm{tr}(zx)\mathcal{L}_{\nu}f(z)\\
& & +\nu\mathcal{L}_{\nu}(\mathrm{tr}(x\nabla
)f)(z) +\mathcal{L}_{\nu}(\mathrm{tr}(s\nabla x\nabla)f)(z)\\
&  =&\nu\mathrm{tr}(x)\mathcal{L}_{\nu}f(z) +\mathcal{L}_{\nu}(\mathrm{tr}%
(sx)f)(z)\\
&& +\mathcal{L}_{\nu}(\mathrm{tr}((sx+xs)\nabla)f)(z) \\
&  &+\nu\mathcal{L}_{\nu}(\mathrm{tr}(x\nabla)f)(z) +\mathcal{L}_{\nu}%
(\mathrm{tr}(s\nabla x\nabla)f)(z) \\
&  =&\mathcal{L}_{\nu}(\mathrm{tr}((\nu x+sx)+(\nu x+sx+xs)\nabla+s\nabla
x\nabla)f)(z)
\end{eqnarray*}
This completes the proof.

\section{Recursion Relations}

We will now state the main recursion relations implied by these
actions. The notation and definition for the Laguerre functions
are in \cite{doz2}, section 7.
\begin{thm}\label{thm5.1}
The Laguerre functions are related by the following recursion relations:

\begin{enumerate}

\item $\mathrm{tr}(-s\nabla\nabla-\nu \nabla+s)\ell_{\mathbf{m}}^{\nu}=(r\nu+2|\mathbf{m}|)\ell_{\mathbf{m}}^{\nu}.$

\item $\frac{1}{2}\mathrm{tr}(s\nabla\nabla+(\nu I+2s)\nabla+(\nu I+s)
)\ell_{\mathbf{m}}^{\nu}(s)=-\sum_{j=1}^{r}%
\begin{pmatrix}
\mathbf{m}\\
\mathbf{m-\gamma_{j}}%
\end{pmatrix}
(m_{j}-1+\nu-(j-1))\ell_{\mathbf{m-\gamma_{j}}}^{\nu} $

\item $\frac{1}{2}\mathrm{tr}(-s\nabla\nabla+(-\nu I+2s)\nabla+(\nu I-s)
)\ell_{\mathbf{m}}^{\nu}=\sum_{j=1}^{r}c_{\mathbf{m}}(j)\ell_{\mathbf{m+\gamma
_{j}}}^{\nu}.$

\end{enumerate}
\end{thm}

\begin{proof}
Let $\xi=\begin{pmatrix}
  0 & 1 \\
  1 & 0
\end{pmatrix}.$  Then by lemma 5.5 of \cite{doz2}
$$\pi_\nu(\xi)q_{\mathbf{m},\nu}=(r\nu +2 \vert
\mathbf{m}\vert)q_{\mathbf{m},\nu},$$ where, by theorem 7.8 of
\cite{doz2},
$\mathcal{L}_\nu(\ell_{\mathbf{m}}^\nu)=\Gamma_\Omega(\mathbf{m}+\nu)q_{\mathbf{m},\nu}.$
From this it follows that
\begin{eqnarray*}
\lambda_\nu(\xi)\ell_{\mathbf{m}}^\nu &=&
\mathcal{L}_\nu^{-1}\pi_\nu(\xi)\mathcal{L}_\nu
\ell_{\mathbf{m}}^\nu\\
&=&\mathcal{L}_\nu^{-1}\pi_\nu(\xi)\Gamma_\nu(\mathbf{m}+\nu)q_{\mathbf{m},\nu}\\
&=& (r\nu+2\vert
\mathbf{m}\vert)\mathcal{L}_\nu^{-1}(\Gamma_\Omega(\mathbf{m}+\nu)q_{\mathbf{m},\nu})\\
&=&(r\nu+2\vert \mathbf{m}\vert)\ell_{\mathbf{m}}^\nu.
\end{eqnarray*}
On the other hand, if $a=0$ and $b=1$ in Theorem \ref{th2.1}, then
$$ \lambda_\nu(\xi)\ell_{\mathbf{m}}^\nu =\tr(-s\nabla\nabla
-\nu\nabla +s)\ell_\mathbf{m}^\nu$$ and part (1) follows.

The vector $\xi$ is in the center of $\mathfrak{k}_\mathbb{C}$ and
$\operatorname{ad}(\xi)$ acts by $-2I$ on $\mathfrak{p}^+$ and
$2I$ on $\mathfrak{p}^-$, where $I$ is the identity operator. Let
$$L^2_k(\Omega,d\mu_\nu)=\set{f\in
L^2(\Omega,d\mu_\nu):\lambda_\nu(\xi)f=(r\nu+2k)f}.$$ Since
$\lambda_\nu$ is an irreducible highest weight representation it
is well known that $L^2_k(\Omega,d\mu_\nu)$ is finite dimensional,
nonzero if $k\ge0$, and $L^2(\Omega,d\mu_\nu)=\oplus
L^2_k(\Omega,d\mu_\nu).$ Furthermore, part (1) implies that
$\ell_\mathbf{m}^\nu \in
L^2_{\vert\mathbf{m}\vert}(\Omega,d\mu_\nu),$ for all
$\mathbf{m}\in \Lambda$ (c.f. (2.15) of \cite{doz2}). For $x\in
\mathfrak{p}^+$ and $f\in L^2_k(\Omega,d\mu_\nu)$ we have
\begin{eqnarray*}
\lambda_\nu(\xi)\lambda_\nu(x)f&=& \lambda_\nu(x)\lambda_\nu(\xi)f
+\lambda_\nu([\xi,x])f\\
&=&(r\nu +2k)\lambda_\nu(x)f -2 \lambda_\nu(x)f\\
&=&(r\nu +2(k-1))\lambda_\nu(x)f.
\end{eqnarray*}
This implies that $\lambda_\nu(x)f\in L^2_{k-1}(\Omega,d\mu_\nu)$.
Similarly, for $x\in \mathfrak{p}^-$, we have $\lambda_\nu(x)f \in
L^2_{k+1}(\Omega,d\mu_\nu).$

Now let $X^+=\begin{pmatrix}
  -1 & -1 \\
  1 & 1
\end{pmatrix}\in \mathfrak{p}^+$, $X^-=\begin{pmatrix}
  1 & -1 \\
  1 & -1
\end{pmatrix} \in \mathfrak{p}^-$, and $Z=\frac{1}{2}(X^--X^+)=\begin{pmatrix}
  1 & 0 \\
  0 & -1
\end{pmatrix}.$ Then by theorem 7.9 of \cite{doz2} and its proof
$$\lambda_\nu(Z)\ell_{\mathbf{m}}^\nu =-\sum
_{j=1}^{r}{\binom{m}{m-\gamma_{j}}}(m_{j}-1+\nu-(j-1))\ell_{m-\gamma_{j}%
}^{\nu}+\sum_{j=1}^{r}c_{m}(j)\ell_{m+\gamma_{j}}^{\nu}.$$ If
$P_k$ denotes the orthogonal projection of $L^2(\Omega,d\mu_\nu)$
onto $L^2_k(\Omega,d\mu_\nu)$ then $$ \frac
{-1}{2}\lambda_\nu(X^+)\ell_\mathbf{m}^\nu
=P_{\vert\mathbf{m}\vert-1}\lambda_\nu(Z)\ell_\mathbf{m}^\nu=-\sum
_{j=1}^{r}{\binom{m}{m-\gamma_{j}}}(m_{j}-1+\nu-(j-1))\ell_{m-\gamma_{j}%
}^{\nu}$$ and $$ \frac {1}{2}\lambda_\nu(X^-)\ell_\mathbf{m}^\nu
=P_{\vert\mathbf{m}\vert+1}\lambda_\nu(Z)\ell_\mathbf{m}^\nu=
\sum_{j=1}^{r}c_{m}(j)\ell_{m+\gamma_{j}}^{\nu}.$$ Formulas (2)
and (3) now follow from theorem \ref{th2.1} part (2) and (3),
respectively, by setting $x=1$.

\end{proof}

\end{document}